\documentclass[12pt,a4paper]{article}
\usepackage[a4paper,margin=15mm]{geometry}
\usepackage[T1]{fontenc}
\usepackage{amsmath,amssymb,amsthm,booktabs}
\usepackage[hidelinks]{hyperref}
\newtheorem{theorem}{Theorem}[section]
\newtheorem{lemma}[theorem]{Lemma}

\newtheorem{corollary}[theorem]{Corollary}
\newcommand{\coop}{\chi_{\mathrm{cp}}}
\newcommand{\Cell}{\mathcal C_\ell}

\theoremstyle{definition}
\newtheorem{question}[theorem]{Question}

\begin{document}

\title{Cooperative colorings of cycles and uniform cycle factors\thanks{Supported by the National Natural Science Foundation of China (Nos. 12501495, 12501496) and Shaanxi Fundamental Science Research Project for Mathematics and Physics (Grant No.\,25JSQ050).}}
\author{Xuqing Bai \and Zhipeng Gao \and Bi Li\thanks{Corresponding author.} \and Xin Zhang}
\date{}
\maketitle
\begin{center}
School of Mathematics and Statistics, Xidian University, Xi'an 710071, China\\[2pt]
\small\texttt{baixuqing@xidian.edu.cn}, \texttt{gaozhipeng@xidian.edu.cn},\\
\texttt{libi@xidian.edu.cn}, \texttt{xzhang@xidian.edu.cn}
\end{center}
\begin{abstract}
A cooperative coloring of graphs on a common vertex set is a choice of
one independent set in each graph such that their union covers the vertex
set. We give a combinatorial proof that three cycles on a common vertex
set always have a cooperative coloring. We also prove that, for every
$\ell\geq4$ with $\ell\equiv1\pmod3$, there are three graphs on a common
vertex set, each a disjoint union of cycles of length $\ell$, with no
cooperative coloring. Together with known bounds, this determines, for every $\ell\geq3$,
the least number $m$ such that every $m$ graphs on a common vertex set, each
a disjoint union of cycles of length $\ell$, have a cooperative coloring.
For $\ell\geq3$, this number is four when $\ell\equiv1\pmod3$ and three
otherwise. These results answer two questions of Aharoni, Holzman, Howard,
and Spr\"ussel.
\end{abstract}

\section{Introduction}
All graphs in this paper are finite and simple, except for an auxiliary
multigraph used to apply a coloring theorem. A set of vertices is independent
if it contains no adjacent pair. Let $G_1,\ldots,G_m$ be graphs on a common
vertex set $V$. A {\it cooperative coloring} of $G_1,\ldots,G_m$ is a tuple
$(I_1,\ldots,I_m)$ such that $I_i$ is independent in $G_i$ for each
$i\in\{1,\ldots,m\}$ and
\[
V=I_1\cup\cdots\cup I_m.
\]
The graphs need not be distinct. The sets $I_1,\ldots,I_m$ may overlap. 
By deleting vertices from the sets $I_i$ if necessary, we may assume
that $I_1,\ldots,I_m$ form a partition of $V$, that is,
\[
V=I_1\mathbin{\dot\cup}\cdots\mathbin{\dot\cup}I_m.
\]

A cooperative coloring exists if and only if there is a map
$\varphi\colon V\to\{1,\ldots,m\}$ such that, for every
$i\in\{1,\ldots,m\}$,
\begin{equation}\label{eq:map}
uv\in E(G_i)\quad\Longrightarrow\quad
(\varphi(u),\varphi(v))\ne(i,i).
\end{equation}
If $G_1=\cdots=G_m=G$, such a map is a proper coloring of $G$ with at
most $m$ colors. For a graph class $\mathcal C$, the {\it cooperative
chromatic number} $\coop(\mathcal C)$ is the least positive integer $m$
such that every $m$ graphs in $\mathcal C$ on a common vertex set have
a cooperative coloring.


Aharoni, Holzman, Howard, and Spr\"ussel~\cite{AHHS} studied
cooperative colorings using independent systems of representatives,
also called independent transversals. An {\it independent transversal}
of a partition of a graph's vertex set is an independent set containing
exactly one vertex from each part. The independent transversal theorem
of Haxell~\cite{Haxell} implies that, for every positive integer
$\Delta$, any $2\Delta$ graphs on a common vertex set, each with maximum
degree at most $\Delta$, have a cooperative coloring. Aharoni et
al.~\cite{AHHS} showed that, for every $\Delta\geq2$, there are
$\Delta+1$ graphs of maximum degree at most $\Delta$ on a common vertex
set with no cooperative coloring. A result of Loh and
Sudakov~\cite{LohSudakov} implies that $(1+o(1))\Delta$ graphs of
maximum degree at most $\Delta$ suffice as $\Delta\to\infty$.

Cooperative colorings have also been studied for various graph classes.
A result of Aharoni, Berger, and Ziv~\cite{ABZ} implies that any
$\Delta+1$ chordal graphs of maximum degree at most $\Delta$ on a common
vertex set have a cooperative coloring. Aharoni, Berger, Chudnovsky,
Havet, and Jiang~\cite{ABCHJ} obtained bounds on the number of graphs
needed to guarantee a cooperative coloring when the graphs are trees or
bipartite graphs of bounded maximum degree. 
They also observed that, for every \(\Delta\ge2\), the minimum number of graphs needed to guarantee a cooperative coloring is the same for trees and forests of maximum degree at most \(\Delta\).
Bradshaw~\cite{Bradshaw}
improved the lower bound for forests and determined the asymptotic
behavior for star forests. 
Bai, Li, Xu, and Zhang~\cite{BLXZ} constructed three trees of maximum degree at most three on a common vertex set with no cooperative coloring. Together with the preceding upper bound for chordal graphs, this shows that four is the minimum number of trees of maximum degree at most three needed to guarantee a cooperative coloring.
They also obtained results for wheels and for
graphs whose components are balanced complete bipartite graphs or
generalized theta graphs.

Cooperative coloring is closely related to adapted coloring, studied
by Hell and Zhu~\cite{HellZhu}. Kostochka and Zhu~\cite{KostochkaZhu}
introduced and studied its list-coloring extension. Bradshaw and
Masa\v r\'ik~\cite{BM} studied single-conflict colorings, a generalization of adapted colorings, and obtained
an upper bound for cooperative colorings of families of degenerate
graphs with bounded maximum degree.

Cooperative coloring has also been studied for matroids and
hypergraphs. Bartnicki, Czerwi\'nski, Grytczuk, and
Miechowicz~\cite{BCGM} proved that $k$ matroids on the same ground set
have a cooperative coloring whenever each of them is $k$-colorable.
Bai, Li, Liu, and Zhang~\cite{BLLZ} studied degree conditions for
cooperative colorings of uniform hypergraphs and gave a combinatorial
proof that two $k$-uniform tight cycles on a common vertex set have a
cooperative coloring when $k\geq3$.

We now turn to graphs of maximum degree two, focusing on cycles and
disjoint unions of cycles. 
Aharoni, Holzman, Howard, and Spr\"ussel~\cite{AHHS} proved that three
cycles on the same vertex set have a cooperative coloring. Their proof
uses topological methods. They also posed the following two questions.

\begin{question}\label{question:cycles}
Is there a combinatorial proof that any three cycles on the same vertex
set have a cooperative coloring?
\end{question}

\begin{question}\label{question:factors}
Let $\ell\geq3$. Do any three graphs on the same vertex set, each of
which is a disjoint union of cycles of length $\ell$, have a cooperative
coloring?
\end{question}

Haxell and Wdowinski~\cite{HW} later proved that the disjoint union of
three cycles, each of length congruent to one modulo three, has a
partition into triples with no independent transversal.

We answer both questions. For Question~\ref{question:cycles}, we give a
combinatorial proof using the cycle-plus-triangles theorem of Fleischner
and Stiebitz~\cite{FS}. An elementary combinatorial proof of that theorem
was given by Sachs~\cite{Sachs}. For Question~\ref{question:factors},
the results of Aharoni, Holzman, Howard, and Spr\"ussel~\cite{AHHS} imply an affirmative answer
when $\ell\not\equiv1\pmod3$. We prove that the answer is negative when
$\ell\equiv1\pmod3$.

We state these results formally below. The first is the theorem of Aharoni et al.~\cite{AHHS} for which we provide a combinatorial proof.

\begin{theorem}\label{thm:cycles}
Let $C_1,C_2,C_3$ be cycles on the same vertex set $V$. There is a
partition
\[
V=I_1\mathbin{\dot\cup}I_2\mathbin{\dot\cup}I_3
\]
such that $I_i$ is independent in $C_i$ for each $i\in\{1,2,3\}$.
\end{theorem}

For $\ell\geq3$, let $\Cell$ be the class of graphs whose components
are all cycles of length $\ell$. We call these graphs $C_\ell$-factors.

\begin{theorem}\label{thm:factors}
For every integer $\ell\geq4$ with $\ell\equiv1\pmod3$, there are three
$C_\ell$-factors on a common vertex set with no cooperative coloring.
\end{theorem}

For every $\ell\geq3$, two $C_\ell$-factors on a common vertex set
need not have a cooperative coloring. We prove this in
Subsection~\ref{sec3.3}. 
The results of Aharoni, Holzman, Howard, and Spr\"ussel~\cite{AHHS} imply that
$\coop(\Cell)\leq4$ for every $\ell\geq3$, and that
$\coop(\Cell)\leq3$ when $\ell\not\equiv1\pmod3$.
Together with these
upper bounds, the lower bound above and Theorem~\ref{thm:factors} give
the following classification.

\begin{corollary}\label{cor:classification}
For every integer $\ell\geq3$,
\[
\coop(\Cell)=
\begin{cases}
4,&\ell\equiv1\pmod3,\\
3,&\ell\not\equiv1\pmod3.
\end{cases}
\]
\end{corollary}

The rest of the paper is organized as follows. Section~\ref{sec:cycles}
gives a combinatorial proof of Theorem~\ref{thm:cycles}.
Section~\ref{sec:factors} proves Theorem~\ref{thm:factors} and
Corollary~\ref{cor:classification}.

\section{Three cycles on a common vertex set}\label{sec:cycles}

We use the cycle-plus-triangles theorem of Fleischner and
Stiebitz~\cite{FS}, for which Sachs~\cite{Sachs} gave an elementary
combinatorial proof.

\begin{theorem}\label{thm:CPT}
Let $H$ be a cycle on $3r$ vertices, where $r\geq1$. Partition its
vertices into triples and add a triangle on each triple. The resulting
graph has a proper three-coloring.
\end{theorem}

An edge of an added triangle may already be an edge of $H$. In the
decomposition form of the theorem, these are treated as parallel edges.
Since parallel edges impose the same coloring condition, the assertion
above also holds for the resulting simple graph.

We use Theorem~\ref{thm:CPT} to obtain an independent transversal in a
linear forest. This will allow us to choose the colors of all but one
vertex in the proof of Theorem~\ref{thm:cycles}. A {\it linear forest}
is a graph whose components are paths, with isolated vertices allowed
as paths.

\begin{lemma}\label{lem:forest}
Let $F$ be a linear forest whose vertices are partitioned into nonempty
parts $B_1,\ldots,B_r$. Suppose every part has size three, except that
at most two parts may have size two. Then these parts have an independent
transversal in $F$.
\end{lemma}

\begin{proof}
Add one isolated vertex to each part of size two, and denote the
enlarged parts by $Q_1,\ldots,Q_r$. Every enlarged part has size three.
At most two new vertices have been added, and the enlarged graph is
still a linear forest.

List the vertices of each path in path order, and concatenate the lists.
Join consecutive vertices in the resulting list and join its last vertex
to its first. Since there are $3r\geq3$ vertices, these edges form a
cycle $H$ containing every edge of $F$. Add a triangle on each $Q_i$
and apply Theorem~\ref{thm:CPT} to obtain a proper three-coloring.

At most two of the colors occur on the new vertices. Choose a color
absent from them, and let $S$ be its color class. Every triangle uses
all three colors, so $S$ contains exactly one vertex of each $Q_i$.
This vertex lies in $B_i$, since it is not new. The set $S$ is
independent in $H$, and hence in $F$. Therefore $S$ is the required
transversal.
\end{proof}

\begin{proof}[Proof of Theorem~\ref{thm:cycles}]
Choose $v\in V$, and let $u,w$ be its two neighbors in $C_1$. Make
disjoint copies of $C_1-\{u,v,w\}$, $C_2-v$, and $C_3-v$, with vertex
sets
\[
W_1=\{(x,1)\mid x\in V\setminus\{u,v,w\}\},\qquad
W_i=\{(x,i)\mid x\in V\setminus\{v\}\}\quad(i=2,3).
\]
Adjacency in each copy is inherited from its cycle. Let $F$ be the
disjoint union of the copies. Deleting one vertex from a cycle leaves
a path. Also, $u,v,w$ are consecutive in $C_1$, so deleting them leaves
a path, possibly empty. Thus $F$ is a linear forest.

For every $x\in V\setminus\{v\}$, let
\[
B_x=\{(x,i)\mid(x,i)\in V(F)\}.
\]
These sets partition $V(F)$. The two parts $B_u,B_w$ have size two,
while every other part has size three. By Lemma~\ref{lem:forest}, they
have an independent transversal $S$.

Set $\varphi(v)=1$. For each $x\ne v$, set $\varphi(x)=i$ when $(x,i)$
is the unique vertex in $S\cap B_x$. If $xy\in E(C_i)$ with $x,y\ne v$
and both endpoints had color $i$, then $(x,i)$ and $(y,i)$ would be
adjacent members of $S$, contrary to its independence.

An edge incident with $v$ in $C_2$ or $C_3$ cannot have both endpoints
colored with its cycle index, since $\varphi(v)=1$. The two edges
incident with $v$ in $C_1$ have other endpoints $u,w$. Neither can
receive color one because its copy in $W_1$ was removed. Therefore
$\varphi$ satisfies~\eqref{eq:map}. Its three color classes give the
required partition.
\end{proof}

\section{\texorpdfstring{$C_\ell$}{C-ell}-factors on a common vertex set}\label{sec:factors}

\subsection{The case \texorpdfstring{$\ell\geq7$}{ell >= 7}} 

Throughout this subsection, let $\ell\geq7$ and
$\ell\equiv1\pmod3$. We use the following result of Haxell and
Wdowinski~\cite{HW}.

\begin{lemma}\label{thm:HW}
Let $\ell_1,\ell_2,\ell_3\geq4$ with
$\ell_i\equiv1\pmod3$ for each $i\in\{1,2,3\}$. The disjoint union of
cycles of lengths $\ell_1,\ell_2,\ell_3$ has a partition into triples
with no independent transversal.
\end{lemma}

The following lemma gives a property used in the proof of
Theorem~\ref{thm:factors}. We state only the edges needed for the proof.
Additional edges do not affect the conclusion.

\begin{lemma}\label{lem:seven}
Let $G_1,G_2,G_3$ be graphs on a common vertex set $V$. Suppose that they
contain seven distinct vertices $u,v,w,a,b,c,d$ and the edges
\[
\begin{aligned}
db,bc,ca&\in E(G_1),\\
ab,bu,uv,vw,wc,cd&\in E(G_2),\\
au,uw,wd,bv,vc&\in E(G_3).
\end{aligned}
\]
If a map $\varphi\colon V\to\{1,2,3\}$ satisfies~\eqref{eq:map}, then
at least one of $u,v,w$ has color one.
\end{lemma}

\begin{proof}
Suppose that all of $u,v,w$ have colors in $\{2,3\}$. If
$\varphi(v)=2$, the edges $uv,vw$ of $G_2$ force
$\varphi(u)=\varphi(w)=3$. This contradicts the edge $uw$ of $G_3$.
Hence $\varphi(v)=3$.

The edges $bv,vc$ of $G_3$ give
$\varphi(b),\varphi(c)\in\{1,2\}$. The edge $uw$ of $G_3$ also shows
that at least one of $u,w$ has color two.

First suppose that $\varphi(u)=2$. The edge $bu$ of $G_2$ forces
$\varphi(b)=1$. The edge $bc$ of $G_1$ then forces $\varphi(c)=2$.
The edge $cw$ of $G_2$ forces $\varphi(w)=3$. Now $cd\in E(G_2)$
excludes color two from $d$, and $wd\in E(G_3)$ excludes color three.
Thus $\varphi(d)=1$, contrary to $db\in E(G_1)$.

It remains to consider $\varphi(w)=2$. The edge $cw$ of $G_2$ gives
$\varphi(c)=1$. The edge $bc$ of $G_1$ gives $\varphi(b)=2$. The edge
$bu$ of $G_2$ then gives $\varphi(u)=3$. The edges $ab\in E(G_2)$ and
$au\in E(G_3)$ exclude colors two and three from $a$. Thus
$\varphi(a)=1$, contrary to $ac\in E(G_1)$. Both possibilities are
impossible.
\end{proof}

We write $(x_1,\ldots,x_s)$ for the cycle with edges between consecutive
vertices and between $x_s$ and $x_1$.

\begin{proof}[Proof of Theorem~\ref{thm:factors} for $\ell\geq7$]
By Lemma~\ref{thm:HW}, the disjoint union $H$ of three $\ell$-cycles
has a partition $\mathcal T$ into triples with no independent
transversal. There are $\ell$ triples because $|V(H)|=3\ell$.

For every $T\in\mathcal T$, choose an ordering $T=\{u,v,w\}$ and
introduce the new vertices
\[
a_T,b_T,c_T,d_T,\quad p_{T,1},\ldots,p_{T,\ell-7},
\quad f_{T,1},f_{T,2},f_{T,3}.
\]
All new vertices are outside $V(H)$, and vertices introduced for
different triples are distinct. For a fixed $T$, omit the subscript $T$
when writing its cycles. The sequence of $p$-vertices is empty when
$\ell=7$.

Let $G_1$ contain the three cycles of $H$ and, for each triple $T$,
the cycle
\begin{equation}\label{eq:G1}
(d,b,c,a,p_1,\ldots,p_{\ell-7},f_1,f_2,f_3).
\end{equation}
For the same triple, put the cycles
\begin{equation}\label{eq:G23}
\begin{aligned}
(a,b,u,v,w,c,d,p_1,\ldots,p_{\ell-7})
    &\quad\text{in }G_2,\\
(a,u,w,d,b,v,c,p_1,\ldots,p_{\ell-7})
    &\quad\text{in }G_3.
\end{aligned}
\end{equation}
Each displayed cycle has length $\ell$. The only vertices not yet used
in $G_2$ and $G_3$ are in
\[
F=\bigcup_{T\in\mathcal T}\{f_{T,1},f_{T,2},f_{T,3}\}.
\]
Since $|F|=3\ell$, partition $F$ into three sets of size $\ell$ and
arrange each set as an $\ell$-cycle in both $G_2$ and $G_3$.

Every original vertex lies in exactly one triple. Thus it lies in
exactly one displayed cycle of each of $G_2,G_3$, and in one cycle of
$H$ in $G_1$. Every new vertex also lies in exactly one cycle in each
graph. Therefore $G_1,G_2,G_3$ are $C_\ell$-factors on a common vertex
set. There are $\ell$ new vertices for each of the $\ell$ triples.
Their common order is
\[
3\ell+\ell\cdot\ell=\ell(\ell+3).
\]

Suppose that these factors have a cooperative coloring. By
\eqref{eq:map}, there is a corresponding map $\varphi$. For each
$T=\{u,v,w\}$, the cycles in~\eqref{eq:G1} and~\eqref{eq:G23} contain
all edges required by Lemma~\ref{lem:seven}. Hence $\varphi^{-1}(1)$
meets every triple.

The set $\varphi^{-1}(1)\cap V(H)$ is independent in $H$, because $H$
is a subgraph of $G_1$. Choosing one of its vertices from each triple
gives an independent transversal of $\mathcal T$, contrary to the
choice of $H$ and $\mathcal T$.
\end{proof}

\subsection{The case \texorpdfstring{$\ell=4$}{ell = 4}}

For this case, we use three graphs on 18 vertices. Their cooperative
colorings satisfy a property that will allow us to obtain the required
$C_4$-factors from two copies of these graphs.

Let $W=\{1,\ldots,18\}$, and define $H_1,H_2,H_3$ by the components in
Table~\ref{tab:components}. A bipartition of a graph is a partition into
two independent sets. On a four-cycle, these two sets are the pairs of
opposite vertices. Every independent set of the cycle is contained in
one such pair.

\begin{table}[htbp]
\centering
\renewcommand{\arraystretch}{1.18}
\begin{tabular}{ccc}
\toprule
$H_1$&$H_2$&$H_3$\\
\midrule
$(3,5,4,6)$&$(15,3,16,9)$&$(1,5,14,18)$\\
$(7,9,8,10)$&$(1,6,5,17)$&$(11,4,13,17)$\\
$(11,13,12,14)$&$(2,10,12,11)$&$(3,7,6,8)$\\
$(15,17,16,18)$&$(4,14,13,18)$&$(2,9,10,12)$\\
$\{1\},\{2\}$&$\{7,8\}$&$\{15,16\}$\\
\bottomrule
\end{tabular}
\caption{Components of $H_1,H_2,H_3$. Four-term sequences denote
cycles. Two-element sets denote single edges, and singletons denote
isolated vertices.}
\label{tab:components}
\end{table}

\begin{lemma}\label{lem:eighteen}
Every cooperative coloring $(I_1,I_2,I_3)$ of $(H_1,H_2,H_3)$ satisfies
\[
I_1\cap\{1,2\}\ne\varnothing.
\]
\end{lemma}

\begin{proof}
Suppose instead that $I_1\cap\{1,2\}=\varnothing$. On each four-cycle
component of $H_i$, enlarge the selection in $I_i$ to a pair of
opposite vertices containing it. Keep its selection on isolated vertices
and single-edge components unchanged, and call the resulting independent
set $J_i$. The sets still cover $W$, and
$J_1\cap\{1,2\}=\varnothing$. On each four-cycle, opposite vertices now
either both belong to $J_i$ or both lie outside $J_i$.

We first prove that $15,16\in J_1$. Suppose $15\in J_2$. Then
$16\in J_2$ and $3,9\notin J_2$, by the first cycle of $H_2$. If
$7\in J_1$, then $8\in J_1$ and $9,10\notin J_1$. Covering 9 forces
$9\in J_3$, which excludes 2 and 10 from $J_3$. Since $2\notin J_1$,
covering 2 then forces $2\in J_2$. This excludes 10 from $J_2$, leaving
10 uncovered. Therefore $7,8\notin J_1$.

The vertices 7 and 8 are adjacent in $H_2$, so at least one must belong
to $J_3$. They are opposite in the third cycle of $H_3$, so both belong
to $J_3$. Hence $3,6\notin J_3$. Covering 3 forces $3\in J_1$, which
excludes 6 from $J_1$. Covering 6 now forces $6\in J_2$, which excludes
1 from $J_2$. Since $1\notin J_1$, we obtain $1\in J_3$. But 5 is
adjacent to 3 in $H_1$, to 6 in $H_2$, and to 1 in $H_3$. It is
therefore uncovered, a contradiction.

Thus $15\notin J_2$, and its opposite vertex 16 also lies outside $J_2$.
They are adjacent in $H_3$, so at least one belongs to $J_1$. They are
opposite in $H_1$, giving $15,16\in J_1$. Consequently
$17,18\notin J_1$.

We next prove that $11,13\notin J_3$. If $11\in J_3$, then
$13\in J_3$ and $4,17\notin J_3$. Since $17\notin J_1$, covering 17
gives $17\in J_2$. Its edge to 1 in $H_2$ excludes 1 from $J_2$, so
$1\in J_3$. The edge from 1 to 18 in $H_3$ excludes 18 from $J_3$.
Since $18\notin J_1$, we get $18\in J_2$. Its edge to 4 in $H_2$
excludes 4 from $J_2$, and hence $4\in J_1$. But 5 is adjacent to 4 in
$H_1$, to 17 in $H_2$, and to 1 in $H_3$. Again 5 is uncovered. It
follows that $11\notin J_3$, and the same holds for its opposite vertex
13.

If $11\notin J_1$, then $12\notin J_1$ because these vertices are
opposite in $H_1$. Covering 11 gives $11\in J_2$, which excludes 2 and
12 from $J_2$. Since $2\notin J_1$, it follows that $2\in J_3$, which
excludes 12 from $J_3$. This leaves 12 uncovered. Therefore
$11,12\in J_1$, and $13,14\notin J_1$.

Since $13\notin J_3$, covering 13 forces $13\in J_2$. This excludes 14
from $J_2$, so $14\in J_3$. The vertex 18 is adjacent to 15 in $H_1$,
to 13 in $H_2$, and to 14 in $H_3$. It is therefore uncovered, a final
contradiction. The enlargement did not change membership of 1 and 2 in
the first set, so the assertion holds for $I_1$.
\end{proof}

\begin{proof}[Proof of Theorem~\ref{thm:factors} for $\ell=4$]
Take a second copy $H_1',H_2',H_3'$ on the disjoint vertex set
$W'=\{1',\ldots,18'\}$. Add the cycle $(1,1',2,2')$ to
$H_1\cup H_1'$ to obtain $G_1$. In $H_2\cup H_2'$, complete the two
single-edge components to the cycle $(7,8,7',8')$, obtaining $G_2$. In
$H_3\cup H_3'$, complete the two single-edge components to the cycle
$(15,16,15',16')$, obtaining $G_3$. All other components remain
unchanged. Every vertex now lies in one four-cycle of each graph, so
these are three $C_4$-factors on 36 vertices.

Suppose $(I_1,I_2,I_3)$ is a cooperative coloring of these factors.
Restricting it to $W$ gives a cooperative coloring of
$(H_1,H_2,H_3)$, since all edges of those graphs remain present.
Restriction to $W'$ gives a cooperative coloring of the second copy.
Lemma~\ref{lem:eighteen} therefore gives
\[
I_1\cap\{1,2\}\ne\varnothing,\qquad
I_1\cap\{1',2'\}\ne\varnothing.
\]
The two displayed pairs are the bipartition classes of the added
four-cycle in $G_1$. Every vertex of either pair is adjacent to every
vertex of the other pair. An independent set of $G_1$ cannot meet both,
a contradiction.
\end{proof}

\subsection{The cooperative chromatic number} \label{sec3.3}

To prove Corollary~\ref{cor:classification}, we first establish a lower
bound of three.

\begin{lemma}\label{lem:two}
For every integer $\ell\geq3$, there are two $C_\ell$-factors on a
common vertex set with no cooperative coloring.
\end{lemma}

\begin{proof}
For $\ell=3$, take two copies of a triangle. Each independent set has
size at most one, so two such sets cannot cover its three vertices.

For $\ell\geq4$, consider the paths $P=v_1v_2v_3v_4$ and
$Q=v_1v_3v_2v_4$, which also appear in~\cite{ABCHJ}. We check directly
that they have no cooperative coloring. If $v_2$ receives color one,
the edges $v_1v_2,v_2v_3$ of $P$ force $v_1,v_3$ to receive color two.
This contradicts $v_1v_3\in E(Q)$. Thus $v_2$ receives color two. The
edges $v_2v_3,v_2v_4$ of $Q$ then force $v_3,v_4$ to receive color
one, contrary to $v_3v_4\in E(P)$.

On $\{v_1,\ldots,v_\ell\}$, take the two cycles
\[
(v_1,v_2,v_3,v_4,v_5,\ldots,v_\ell),\qquad
(v_1,v_3,v_2,v_4,v_5,\ldots,v_\ell).
\]
For $\ell=4$, each sequence ends at $v_4$. These cycles contain $P$ and
$Q$, respectively. A cooperative coloring of the cycles would restrict
to a cooperative coloring of the two paths, which is impossible. Each
cycle alone is a $C_\ell$-factor.
\end{proof}

\begin{proof}[Proof of Corollary~\ref{cor:classification}]
By Lemma~\ref{lem:two}, we have $\coop(\Cell)\geq3$ for every
$\ell\geq3$.
If $\ell\not\equiv1\pmod3$, the known upper bound
$\coop(\Cell)\leq3$ gives $\coop(\Cell)=3$.
If $\ell\equiv1\pmod3$, then $\ell\geq4$.
Theorem~\ref{thm:factors} gives $\coop(\Cell)\geq4$, while the general
degree bound gives $\coop(\Cell)\leq4$. Hence $\coop(\Cell)=4$.
\end{proof}

\end{document}